\documentclass[12pt]{article}
\usepackage{amsfonts}

\newtheorem{theorem}{Theorem}[section]
\newtheorem{prop}[theorem]{Proposition}
\newtheorem{unnumber}{}

\newtheorem{prepf}[unnumber]{Proof}
\newenvironment{pf}{\prepf\rm}{\endprepf}
\newcommand{\qed}{\qquad$\Box$}

\newcommand{\ee}{\mathord{\mathrm{e}}}
\newcommand{\dd}{\mathord{\mathrm{d}}}

\renewcommand{\wr}{\mathbin{\mathrm{wr}}}
\def\ox{\overline{x}}
\def\oy{\overline{y}}
\def\oR{\overline{R}}
\def\bp{B^\prime}

\begin{document}

\title{Asymptotics for incidence matrix classes}
\author{Peter Cameron, Thomas Prellberg and Dudley Stark\\
\small{School of Mathematical Sciences}\\[-0.8ex]
\small{Queen Mary, University of London}\\[-0.8ex]
\small{Mile End Road, London, E1 4NS  U.K.}\\[-0.8ex]
\small\texttt{$\{$p.j.cameron, t.prellberg, d.s.stark$\}$@qmul.ac.uk}}
\date{
}
\maketitle

\begin{abstract}
We define {\em incidence matrices} to be zero-one matrices
with no zero rows or columns. We are interested in counting
incidence matrices with a given number of ones, irrespective
of the number of rows or columns.
A classification of incidence matrices is considered
for which conditions of symmetry by transposition, 
having no repeated rows/columns,
or identification by permutation of rows/columns are imposed.
We find asymptotics and relationships for the number of matrices
with $n$ ones in some of these classes as $n\to\infty$.
\end{abstract}


\section{Introduction}

In this paper we address the problem: \textit{How many zero-one matrices
are there with exactly $n$ ones?} Note that we do not specify in advance
the number of rows or columns of the matrices. In order to make the answer 
finite, we assume that no row or column of such a matrix consists entirely 
of zeros. We call such a matrix an {\em incidence matrix}.

Rather than a single problem, there are many different problems here,
depending on what symmetries and constraints are permitted. In general,
we define $F_{ijkl}(n)$ to be the number of zero-one matrices with $n$ ones
and no zero rows or columns, subject to the conditions
\begin{itemize}
\item $i=0$ if matrices differing only by a row permutation are identified,
and $i=1$ if not;
\item $j=0$ if matrices with two equal rows are forbidden, and $j=1$ if not;
\item $k=0$ if matrices differing only by a column permutation are identified,
and $k=1$ if not;
\item $l=0$ if matrices with two equal columns are forbidden, and $l=1$ if not.
\end{itemize}
The notation is chosen so that $F_{ijkl}(n)$ is a monotonic increasing
function of each of the arguments $i,j,k,l$.

By transposition, it is clear that $F_{klij}(n)=F_{ijkl}(n)$ for all $i,j,k,l,
n$. So, of the sixteen different functions defined above, only ten are
distinct. However, among the problems with $k=i$ and $l=j$, we may decide
that matrices which are transposes of each other are identified, leading
to four further counting problems $\Phi_{ij}(n)$, for $i,j\in\{0,1\}$.

For example, there are four matrices with $n=2$, as shown:
\[\pmatrix{1&1\cr},\quad \pmatrix{1\cr1\cr},\quad
\pmatrix{1&0\cr0&1\cr},\quad \pmatrix{0&1\cr1&0\cr}.\]
The first has repeated columns and the second has repeated rows. The third
and fourth are equivalent under row permutations or column permutations,
while the first and second are equivalent under transposition.

Table~\ref{values} gives some values of these functions. The values of
$F_{1111}(n)$ are taken from the \emph{On-Line Encyclopedia of Integer 
Sequences}~\cite{oeis}, where this appears as sequence A101370
and $F_{0101}(n)$ appears as sequence A049311,
while the
values of $F_{0011}(n)$ and $F_{0111}(n)$ are obtained from a formula in
Corollary~3.3 in~\cite{klazar} using MAPLE. Other computations 
were done with \textsf{GAP}~\cite{gap}.

\begin{table}[htbp]
\caption{\label{values}Some values of the fourteen functions}
\[
\begin{array}{||c||r|r|r|r|r|r|r|r|r||}
\hline
n & 1 & 2 & 3 & 4 & 5 & 6 & 7 & 8 & 9 \\
\hline
F_{0000}(n)  & 1 & 1 &  2 &   4 &    7 &    16 &&& \\
F_{0010}(n)  & 1 & 1 &  3 &  11 &   40 &   174 &&& \\
F_{1010}(n)  & 1 & 2 & 10 &  72 &  624 &  6522 &&& \\
F_{0001}(n)  & 1 & 2 &  4 &   9 &   18 &    44 &&& \\
F_{0011}(n)  & 1 & 2 &  7 &  28 &  134 &   729 &   4408 &   29256 &    210710\\
F_{1001}(n)  & 1 & 2 &  6 &  20 &   73 &   315 &&& \\
F_{1011}(n)  & 1 & 3 & 17 & 129 & 1227 & 14123 &&& \\
F_{0101}(n)  & 1 & 3 &  6 &  16 &   34 &    90 &    211 &     558 &      1430\\
F_{0111}(n)  & 1 & 3 & 10 &  41 &  192 &  1025 &   6087 &   39754 &    282241\\
F_{1111}(n)  & 1 & 4 & 24 & 196 & 2016 & 24976 & 361792 & 5997872 & 111969552\\
\Phi_{00}(n) & 1 & 1 &  2 &   3 &    5 &    11 &&& \\
\Phi_{10}(n) & 1 & 2 &  8 &  44 &  340 &  3368 &&& \\
\Phi_{01}(n) & 1 & 2 &  4 &  10 &   20 &    50 &&& \\
\Phi_{11}(n) & 1 & 3 & 15 & 108 & 1045 & 12639 & 181553 & 3001997 &  55999767\\
\hline
\end{array}
\]
\end{table}

The counting problems can be re-interpreted in various ways:

{\bf Counting hypergraphs by weight}
Given a hypergraph on the vertex set $\{x_1,\ldots,x_r\}$, with edges
$E_1,\ldots,E_s$ (each a non-empty set of vertices), the \emph{incidence
matrix} $A=(a_{ij})$ is the matrix with $(i,j)$ entry $1$ if $x_i\in E_j$,
and $0$ otherwise. The \emph{weight} of the hypergraph is the sum of the
cardinalities of the edges. Thus $F_{0101}(n)$ is the number of hypergraphs
of weight~$n$ with no isolated vertices, up to isomorphism; and $F_{1101}(n)$
is the number of (vertex)-labelled hypergraphs of weight~$n$. Putting
$k=1$ corresponds to labelling the edges, a less usual notion.
Moreover, putting $l=0$ corresponds to counting \emph{simple} hypergraphs
(those without repeated edges). The condition $j=0$ is less natural in this
respect, but corresponds to forbidding ``repeated vertices'' (pairs of
vertices which lie in the same edges).

{\bf Counting bipartite graphs by edges}
Given a zero-one matrix $A=(A_{ij})$, there is a (simple) bipartite graph
whose vertices are indexed by the rows and columns of $A$, with an edge
from $r_i$ to $c_j$ if $A_{ij}=1$. The graph has a distinguished bipartite
block (consisting of the rows). Thus, $F_{0101}(n)$ and $F_{1111}(n)$ count
unlabelled and labelled bipartite graphs with $n$ edges and a distinguished
bipartite block, respectively (where, in the labelled case, we assume that
the labels of vertices in the distinguished bipartite block come first);
$\Phi_{01}(n)$ counts unlabelled bipartite
graphs with $n$ edges and a distinguished bipartition.

{\bf Counting pairs of partitions, or binary block designs}
A \emph{block design} is a set of \emph{plots} carrying two partitions,
the \emph{treatment partition} and the \emph{block partition}. It is said
to be \emph{binary} if no two distinct points lie in the same part of both
partitions; that is, if the meet of the two partitions is the partition into
singletons. Thus, $F_{0101}(n)$ is the number of binary block designs with
$n$ plots. Putting $i=1$ or $k=1$ (or both) corresponds to labelling treatments
or blocks (or both). Combinatorialists often forbid ``repeated blocks'' (this
corresponds to putting $l=0$) although this is not natural from the point of
view of experimental design. Similarly $j=0$ corresponds to forbidding
``repeated treatments''. The functions $\Phi_{ij}(n)$ count block designs up to
duality (interchanging treatments and blocks), without or with treatment and
block labelling and/or forbidding repeated blocks and treatments.

{\bf Counting orbits of certain permutation groups}
A permutation group $G$ on a set $X$ is \emph{oligomorphic} if the number
$F^*_n(G)$ of orbits of $G$ on $X^n$ is finite for all $n$. Equivalently, 
the number $F_n(G)$ of orbits on ordered $n$-tuples of distinct elements is
finite, and the number $f_n(G)$ of orbits on $n$-element subsets of $X$ is
finite, for all $n$. These numbers satisfy various conditions, including
the following:
  \begin{itemize}
  \item $\displaystyle{F^*_n(G)=\sum_{k=1}^nS(n,k)F_k(G)}$, where $S(n,k)$ are
Stirling numbers of the second kind;
  \item $f_n(G)\le F_n(G)\le n!f_n(G)$, where the right-hand bound is attained
if and only if the group induced on a finite set by its setwise stabiliser
is trivial.
  \end{itemize}

For example, let $S$ be the symmetric group on an infinite set $X$, and $A$
the group of all order-preserving permutations of the rational numbers. Then
$F_n(S)=f_n(S)=f_n(A)=1$ and $F_n(A)=n!$~.

Now if $H$ and $K$ are permutation groups on sets $X$ and $Y$, then the direct
product $H\times K$ acts coordinatewise on the Cartesian product $X\times Y$.
It is easy to see that $F_n^*(H\times K)=F_n^*(H)F_n^*(K)$.

Let $(x_1,y_1)$, \dots, $(x_n,y_n)$ be $n$ distinct elements of $X\times Y$.
If both $X$ and $Y$ are ordered, then the set of $n$ pairs can be described by
a matrix with $n$ ones in these positions, where the rows and columns of the
matrix are indexed by the sets $\{x_1,\ldots,x_n\}$ and $\{y_1,\ldots,y_n\}$
respectively (in the appropriate order). Moreover, if $X$ is not ordered,
then we can represent the set of pairs as the equivalence class of this
matrix under row permutations, and similarly for columns. Thus
\[F_{0101}(n)=f_n(S\times S),\quad F_{1101}(n)=f_n(S\times A),\quad
F_{1111}(n)=f_n(A\times A).\]

Moreover, the wreath product $H\wr C_2$ is the permutation group on $X^2$ 
generated by $H\times H$ together with the permutation 
$\tau:(x_1,x_2)\mapsto(x_2,x_1)$.
The effect of $\tau$ is to transpose the matrix representing an orbit. So
\[\Phi_{01}(n)=f_n(S\wr C_2),\qquad \Phi_{11}(n)=f_n(A\wr C_2).\]

Discussion of this ``product action'' can be found in~\cite{cgm} and~\cite{MM}.

It is not clear how forbidding repeated rows or columns can be included in
this interpretation.

\section{The asymptotics of $F_{1111}(n)$}

We will use both $F(n)$ and 
$F_{1111}(n)$ to denote the number of incidence
matrices with $n$ ones. This is the
largest of our fourteen functions, so its value gives an upper bound for 
all the others. Indeed, we will see later that $F_{ijkl}(n)=o(F_{1111}(n))$
for $(i,j,k,l)\ne(1,1,1,1)$.

It is possible to compute this function explicitly. For fixed~$n$, let
$m_{ij}$ be the number of $i\times j$ matrices with $n$ ones (and no zero
rows or columns). We set $m_{0,0}(0)=1$ and $F(0)=1$. Then
\begin{equation}\label{firstdisp}
\sum_{i\le k}\sum_{j\le l}{k\choose i}{l\choose j}m_{ij}={kl\choose n},
\end{equation}
so by M\"obius inversion,
\begin{equation}\label{mobius}
m_{kl}=\sum_{i\le k}\sum_{j\le l}(-1)^{k+l-i-j}{k\choose i}{l\choose j}
{ij\choose n},
\end{equation}
and then
\begin{equation}\label{fsum}
F_{1111}(n)=\sum_{i\le n}\sum_{j\le n}m_{ij}.
\end{equation}

For sequence $a_n$, $b_n$, we use the notation $a_n\sim b_n$
to mean $\lim_{n\to\infty}a_n/b_n=1$.
It is clear from the argument above that
\[F_{1111}(n)\le{n^2\choose n}\sim\frac{1}{\sqrt{2\pi n}}(n\ee)^n,\]
and of course considering permutation matrices shows that
\[F_{1111}(n)\ge n!\sim\sqrt{2\pi n}\left(\frac{n}{\ee}\right)^n.\]
\begin{theorem}\label{F1111}
\[
F_{1111}(n)\sim
\frac{n!}4\ee^{-\frac12(\log2)^2}\frac1{(\log2)^{2n+2}}\;.
\]
\label{f14th}
\end{theorem}

We remark that for $n=10$, the asymptotic expression is about $2.5\%$
less than the actual value of $2324081728$.

We have three different proofs of 
Theorem \ref{F1111}. One proof will be given in its entirety and the
other two will be briefly sketched. Their full details can be found
in \cite{cps}. We use the method of the 
first proof to bound $F_{1101}$ in Section~5.
The ideas behind the third proof
lead to a random algorithm for generating incidence matrices counted by 
$F_{1111}(n)$ and by
$\Phi_{11}(n)$. The random algorithm provides an independent
proof of the expression for $F_{1111}(n)$ used in the first proof.\\

\paragraph{First proof}
This proof uses a procedure which, when successful,
generates an incidence matrix uniformly
at random from all incidence matrices. The probability of success can
be estimated and the asymptotic formula for $F_{1111}(n)$ results.

Let $R$ be a binary relation on a set $X$. We say $R$ is {\em reflexive}
if $(x,x)\in R$ for all $x\in X$. We say $R$ is {\em transitive}
if $(x,y)\in R$ and $(y,z)\in R$ implies $(x,z)\in R$.
A {\em partial preorder} is a relation $R$ on $X$ which is reflexive
and transitive. A relation $R$ is said to satisfy {\em trichotomy}
if, for any $x,y\in X$, one of the cases $(x,y)\in R$, $x=y$, or $(y,x)\in R$
holds. We say that $R$ is a {\em preorder} if it is a partial preorder
that satisfies trichotomy. The members of $X$ are said to be the
{\em elements} of the preorder.

A relation $R$ is {\em antisymmetric} if, whenever $(x,y)\in R$ and
$(y,x)\in R$ both hold, then $x=y$. A relation $R$ on $X$ is a
{\em partial order} if it is reflexive, transitive, and antisymmetric.
A relation is a {\em total order}, if it is a partial order
which satisfies trichotomy.
Given a partial preorder $R$ on $X$,
define a new relation $S$ on $X$ by the rule that
$(x,y)\in S$ if and only if both $(x,y)$ and $(y,x)$ belong to $R$. Then
$S$ is an equivalence relation. Moreover, $R$ induces a partial order
$\oR$ on the set of equivalence classes of $S$ in a natural way:
if $(x,y)\in R$, then $(\ox,\oy)\in\oR$, where $\ox$ is the
$S$-equivalence class containing $x$ and similarly for $y$.
We will call an $S$-equivalence class a {\em block}.
If $R$ is a preorder,
then the relation $\oR$ on the equivalence classes of $S$ is a total order.
See Section 3.8 and question 19 of Section 3.13 in \cite{C}
for more on the above definitions and results.
Random preorders are considered in \cite{cs}.

Given a preorder on elements $[n]:=\{1,2,\ldots,n\}$
with $K$ blocks, let $B_1,B_2,\ldots,B_K$
denote the blocks of the preorder. Generate two random preorders uniformly
at random $B_1,B_2,\ldots,B_K$ and $\bp_1,\bp_2,\ldots,\bp_L$.
For each $1\leq i<j\leq n$, define the event $D_{i,j}$ to be
\[
D_{i,j}=\{{\rm for \ each \ of \ the \ two \ preorders \ }
i {\rm \ and \ } j {\rm \ are \ in \ the \ same \ block}\}.
\]
Furthermore, define
\[
W=\sum_{1\leq i<j\leq n}I_{D_{i,j}},
\]
where the indicator random variables are defined by
\[
I_{D_{i,j}}=
\cases{%
1 & if $D_{i,j}$ occurs;\cr
0 & otherwise.\cr}
\]
If $W=0$, then the procedure is successful, in which case
$B_k\cap\bp_l$ consists of either 0 or 1 elements for each
$1\leq k\leq K$ and $1\leq l\leq L$. If the procedure is successful,
then we define
the corresponding $K\times L$ incidence matrix $A$ by
\[
A_{k,l}=
\cases{%
1 & if $B_k\cap \bp_l\neq\emptyset$;\cr
0 & if $B_k\cap \bp_l=\emptyset$.\cr}
\]
It is easy to check that the above definition of $A$ in fact produces
an incidence matrix and that each incidence matrix occurs in $n!$ different 
ways by the construction. It follows that
\[
F_{1111}(n)=\frac{P(n)^2\mathbb{P}(W=0)}{n!},
\]
where $P(n)$ is the number of preorders on $n$ elements if $n\geq 1$
and $P(0)=1$.

It is known (see \cite{Bar}, for example)
that the exponential generating function of $P(n)$
is
\begin{equation}\label{pgenfunc}
\sum_{n=0}^\infty\frac{P(n)}{n!}z^n=\frac{1}{2-\ee^z}.
\end{equation}
The preceding equality implies that
$P(n)$ has asymptotics given by
\begin{equation}\label{Pasymp}
P(n)\sim \frac{n!}{2}\left(\frac{1}{\log2}\right)^{n+1},
\end{equation}
where, given sequences $a_n$, $b_n$ the notation $a_n\sim b_n$ means
that $\lim_{n\to\infty} a_n/b_n=1$.
It remains to find the asymptotics of $\mathbb{P}(W=0)$.

The $r$th falling moment of $W$ is
\begin{eqnarray}
\mathbb{E}(W)_r&=& \mathbb{E} W(W-1)\cdots(W-r+1)\nonumber\nonumber\\
&=&
\mathbb{E}\left(\sum_{{\rm pairs \ }(i_s,j_s){\rm \ different}}
I_{i_1,j_1}\cdots I_{i_r,j_r}\right)\label{falling}\\
&=&
\mathbb{E}\left(\sum_{{\rm all \ }i_s {\rm \ and \ } j_s{\rm \ different}}
I_{i_1,j_1}\cdots I_{i_r,j_r}\right)+
\mathbb{E}\left(\sum\nolimits^\ast
I_{i_1,j_1}\cdots I_{i_r,j_r}\right),\label{falling2}
\end{eqnarray}
with ${\displaystyle \sum\nolimits^\ast}$ defined
to be the sum with all pairs $(i_s,j_s)$ different, but not all
$i_s, j_s$ different.

First we find the asymptotics of the first term in (\ref{falling2}).
For given sequences $i_1,i_2,\ldots,i_r$, $j_1,j_2,\ldots,j_r$,
the expectation
$\mathbb{E}(I_{i_1,j_1}\cdots I_{i_r,j_r})$ is
the number of ways
of forming two preorders on the set of elements
$[n]\setminus\{j_1,j_2,\ldots,j_r\}$ and then for each $s$
adding the element
$j_s$ to the block containing $i_s$ in both preorders
(which ensures that $D_{i_s,j_s}$ occurs for each $s$) and dividing the
result by $P(n)^2$.
Since the number of ways of choosing  
$i_1,i_2,\ldots,i_r$, $j_1,j_2,\ldots,j_r$ equals
$\frac{n!}{2^r(n-2r)!}$,
This gives
\begin{eqnarray*}
\mathbb{E}\left(\sum_{{\rm all \ }i_s {\rm \ and \ } j_s{\rm \ different}}
I_{i_1,j_1}\cdots I_{i_r,j_r}\right)&=&
\frac{n!}{2^r(n-2r)!}
\frac{P(n-r)^2}{P(n)^2}\\
&\sim&
\left(\frac{(\log 2)^2}{2}\right)^r,
\end{eqnarray*}
where we have used (\ref{Pasymp}).

The second term is bounded in the following way. For each sequence
$(i_1,j_1), (i_2,j_2),\ldots, (i_s,j_s)$ in the second term we form
the graph $G$ on vertices $\bigcup_{s=1}^r \{i_s,j_s\}$ with edges
$\bigcup_{s=1}^r \{\{i_s,j_s\}\}$.
Consider the unlabelled graph $G^\prime$ corresponding to $G$
consisting of $v$ vertices and $c$ components.
The number of ways
of labelling $G^\prime$ to form $G$ is bounded by
$n^v$. The number of preorders corresponding
to this labelling is $P(n-v+c)$ because we form a preorder on $n-v+c$ vertices
after which
the vertices in the connected component of $G$ containing a particular
vertex get added to that block.
Therefore, we have
\begin{eqnarray*}
\mathbb{E}\left(\sum\nolimits^\ast
I_{i_1,j_1}\cdots I_{i_r,j_r}\right)&\leq&
\sum_{G^\prime} n^v \frac{P(n-v+c)^2}{P(n)^2}\\
&=&\sum_{G^\prime}O\left(n^{2c-v}\right),
\end{eqnarray*}
where the constant in $O\left(n^{2c-v}\right)$ is uniform over all $G^\prime$
because $v\leq 2r$.
Since at least one vertex is adjacent
to more than one edge, the graph $G$ is not a perfect matching.
Furthermore, each component of $G$ contains at least two vertices.
It follows that $2c<v$ and,
as a result,
\[
\mathbb{E}\left(\sum\nolimits^\ast
I_{i_1,j_1}\cdots I_{i_r,j_r}\right)=O\left(n^{-1}\right).
\]

The preceding analysis shows that
\[
\mathbb{E}(W)_r\sim\left(\frac{(\log 2)^2}{2}\right)^r
\]
for each $r\geq 0$. The method of moments implies that the distribution
converges weakly to the distribution of a Poisson$((\log 2)^2/2)$ distributed
random variable and therefore
\begin{equation}\label{W0}
\mathbb{P}(W=0)\sim\exp\left(-\frac{(\log 2)^2}{2}\right).
\end{equation}
\qed

\paragraph{Second proof} (Sketch)
First, the following
expression for $F_{1111}(n)$
is given in terms of the number of preorders on $k$ elements
as an alternating sum different from and simpler than (\ref{mobius}):
\[F_{1111}(n)=\frac{1}{n!}\sum_{k=1}^ns(n,k)P(k)^2,\]
where and $s(n,k)$ and
$S(n,k)$ are Stirling numbers of the first and second kind respectively.
As in the first proof, the number of pairs of preorders for which the meets
of the blocks form a given $k$-partition of $[n]$ is $k!F(k)$, so
\[P(n)^2=\sum_{k=1}^n S(n,k)k!F(k),\]
and we obtain the result by inversion.
Next, $P(k)$ is replaced by its asymptotic expression (\ref{Pasymp})
with negligible error.
Let
\[F'(n) = \frac{1}{4}\cdot\frac{1}{n!}\sum_{k=1}^n s(n,k)(k!)^2c^{k+1},\]
where $c=1/(\log2)^2$ is as in the statement of the theorem. As we have argued,
$F(n)\sim F'(n)$.

Now, $(-1)^{n-k}s(n,k)$ is the number of permutations in the symmetric group
$S_n$ which have $k$ cycles. So we can write the formula for $F'(n)$ as a sum
over~$S_n$, where the term corresponding to a permutation with $k$ cycles is
$(-1)^{n-k}(k!)^2c^{k+1}$. In particular, the identity permutation gives
us a contribution
\[g(n) = \frac{1}{4}\,n!\,c^{n+1}.\]
To show that $F'(n)\sim Cg(n)$ as $n\to\infty$, where
$C=\exp(-(\log2)^2/2)$,
we write $F'(n)=F'_1(n)+F'_2(n)+F'_3(n)$, where the three
terms are sums over the following permutations:
\begin{description}
\item{$F'_1$:} all involutions (permutations with $\sigma^2=1$);
\item{$F'_2$:} the remaining permutations with $k\ge\lceil n/2\rceil$;
\item{$F'_3$:} the rest of $S_n$.
\end{description}
A further argument shows that 
$F'_1(n)\sim Cg(n)$, while $F'_2(n)=o(g(n))$ and $F'_3(n)=o(g(n))$.
\qed

\paragraph{Third proof} (Sketch)
If one is interested
in asymptotic enumeration of $F(n)$, the formula (\ref{mobius}), being a double
sum over terms of alternating sign, is on first sight rather unsuitable
for an asymptotic analysis.
We present a derivation of the asymptotic form
of $F(n)$ based on the following elegant and elementary identity. (This
identity and equation~(\ref{mobius}) were also derived in~\cite{MM}.)

\begin{prop}
\begin{equation}
\label{ident1}
F(n)=\sum_{k=0}^\infty\sum_{l=0}^\infty\frac1{2^{k+l+2}}{kl\choose n}\;.
\end{equation}
\end{prop}

\begin{pf}
Insert
\begin{equation}\label{sum1}
1=\sum_{k=i}^\infty\frac1{2^{k+1}}{k\choose i}=
\sum_{l=j}^\infty\frac1{2^{l+1}}{l\choose j}
\end{equation} into (\ref{fsum}) and resum using (\ref{firstdisp}).\qed
\end{pf}

The sum in (\ref{ident1}) is dominated by terms where $kl\gg n$.
In this regime, using
$$
{kl\choose n}\sim\frac{(kl)^n}{n!}\ee^{-\frac{n^2}{2kl}}
$$
and approximating the sum in (\ref{ident1}) by an integral
(cf.\ Euler-Maclaurin) leads to
\begin{eqnarray*}
F(n)&\sim&\frac1{4n!}\int \dd k\int \dd l\;\frac{(kl)^n}{2^{k+l}}\ee^{-\frac{n^2}{2kl}}\\
&=&\frac{n^{2n+2}}{4n!}\int \dd\kappa\int \dd\lambda\;
\ee^{n(\log\kappa-\kappa\log2)}\ee^{n(\log\lambda-\lambda\log2)}
\ee^{-\frac1{2\kappa\lambda}}
\end{eqnarray*}
For $n$ large, the integrals are dominated by a small neighborhood around
their respective saddles. As $\ee^{-\frac1{2\kappa\lambda}}$ is independent
of $n$, we can treat the integrals separately.
Using $w(\kappa)=\log\kappa-\kappa\log2$, the saddle
$\kappa_s=\frac1{\log2}$ is determined from $w'(\kappa_s)=0$
($\lambda_s=\frac1{\log2}$ analogously). Approximating the integrals by
a Gaussian around the saddle point gives
\begin{eqnarray*}
F(n)&\sim&\frac{n^{2n+2}}{4n!}
\ee^{nw(\kappa_s)}\sqrt{\frac{2\pi}{n|w''(\kappa_s)|}}
\ee^{nw(\lambda_s)}\sqrt{\frac{2\pi}{n|w''(\lambda_s)|}}
\ee^{-\frac1{2\kappa_s\lambda_s}}\\
&=&\frac{n^{2n+2}}{4n!}
\left(\ee^{n(\log\log2-1)}\sqrt{\frac{2\pi}{n(\log2)^2}}\right)^2
\ee^{-\frac12(\log2)^2}
\end{eqnarray*}
which simplifies to the desired result.
\qed

\section{Generating random incidence matrices}
It is easily shown that (\ref{pgenfunc}) implies that
$$
P(n)=\sum_{k=0}^\infty \frac{k^n}{2^{k+1}}.
$$
Hence, the distribution $\pi_k$ on the natural numbers defined by
$$
\pi_k=\frac{k^n}{P(n)2^{k+1}}
$$
is a probability distribution.
The following way of generating preorders uniformly at random
was given in \cite{MB}.
\begin{theorem}[Maassen, Bezembinder]
Let $A$ be a set of $n$ elements, $n\geq 1$. Let a random preorder $R$
be generated by the following algorithm:
\begin{enumerate}
\item[(i)] Draw an integer-valued random variable $K$ according to the probability distribution $\pi_k$.

\item[(ii)] To each $a\in A$ assign a random score $X_a$ according to the
uniform distribution on $\{1,2,\ldots,K\}$.

\item[(iii)] Put aRb if and only if $X_a\leq X_b$.
\end{enumerate}
Then all of the $P(n)$ possible preorders on $A$ are obtained
with the same probability $P(n)$.
\end{theorem}

Incidence matrices counted
by $F_{1111}(n)$ can be generated uniformly at random by a similar algorithm.
Define a integer valued joint probability distribution function
$\rho_{k,l}$ by
$$
\rho_{k,l}=\frac{1}{F_{1111}(n)}{{kl}\choose n}2^{-k-l-2}.
$$
\begin{theorem}\label{Fgen}
The following algorithm generates a random incidence matrix counted by
$F_{1111}(n)$.
\begin{enumerate}
\item[(i)] Draw integer-valued random variables $K$ and $L$
according to the joint probability distribution $\rho_{k,l}$.

\item[(ii)] Choose a 0-1 matrix with $K$ rows, $L$ columns,
$n$ 1's and $KL-n$ 0's uniformly at random.

\item[(iii)] Delete all rows and columns for which all entries are 0.
\end{enumerate}
\end{theorem}

\begin{pf}
Denote a 0-1 matrix with $k$ rows, $l$ columns, and $n$ 1's a $(k,l)$-matrix.
Denote an incidence matrix with $i$ rows, $j$ columns, and $n$ 1's a 
$(i,j)$-incidence matrix.  
Now, every $(k,l)$-matrix is generated with equal probability 
$$
\frac{\rho_{k,l}}{{kl\choose n}} = \frac{2^{-k-l-2}}{F_{1111}(n)}
$$
and every $(i,j)$-incidence matrix is generated from ${k\choose i}{l\choose j}$
$(k,l)$-matrices.
Averaging over the probability distribution, it follows that every
$(i,j)$-incidence matrix is generated with probability
$$
p(i,j)={kl\choose n}^{-1}\sum_{k,l} {k\choose i}{l\choose j}\rho_{k,l}
$$
Using (\ref{sum1}),
this sum simplifies to $p(i,j)=1/F_{1111}(n)$.
\qed
\end{pf}

\section{Counting symmetric matrices}

In this section we find the asymptotics for $\Phi_{11}(n)$ and show:
\begin{prop}
$\Phi_{11}(n)\sim\frac{1}{2}F_{1111}(n)$.
\end{prop}

\begin{pf}
Clearly we have $\Phi_{11}(n)=\frac{1}{2}(F_{1111}(n)+S_{11}(n))$, where
$S_{11}(n)$ is the number of symmetric matrices with $n$ ones having no zero
rows or columns, where repeated rows or columns are allowed and row or column
permutations are not permitted. So it suffices to show that 
$S_{11}(n)=o(F_{1111}(n))$.

Now let $I(n)$ be the number of solutions of $\sigma^2=1$ in the symmetric
group $S_n$. Then we have
\[ I(n) \le S_{11}(n) \le I(n)P(n)/n! .\]
The lower bound is clear by considering symmetric permutation
matrices. For the upper bound, our analysis of $F_{1111}(n)$ shows that
$n!S_{11}(n)$ is the number of pairs $(R_1,R_2)$ of preorders on
$\{1,\ldots,n\}$ such that no two points $i$ and $j$ lie in the same block
for both preorders, and additionally such that $R_1$ and $R_2$ are
interchanged by some involution $\sigma$ of $\{1,\ldots,n\}$ (corresponding
to transposition of the matrix). 
So instead of choosing $R_1$ and $R_2$, we can choose $R_1$ and $\sigma$ and 
let $R_2=R_1^\sigma$; there are $P(n)I(n)$ choices, 
and this is an overcount because of the extra condition that must hold on 
$(R_1,R_2)$.

Now $I(n)$ is just a little  larger than $\sqrt{n!}$: in fact,
\[I(n) \sim \frac{n^{n/2}}{\sqrt{2}\,\mathrm{e}^{n/2-\sqrt{n}+1/4}}\]
(see \cite[p.~347]{Ben}). We have seen that $P(n)/n!\sim A(1/\log2)^n$.
So the conclusion follows from Theorem~\ref{f14th}.\qed
\end{pf}

It is possible to show that the upper bound for $S_{11}(n)$ is correct, 
apart from a constant factor:

\begin{prop}
$S_{11}(n)\sim C_s\cdot I(n)P(n)/n!$, where 
$C_s={\textstyle{\frac{1}{2}}}\mathrm{e}^{-(\log 2)^2/4}
\approx0.44341.$ In other words, if we choose randomly a preorder $R$ and an
involution $\sigma\in S_n$, the probability that no two points lie in the same
part in both $R$ and $R^\sigma$ tends to $C_s$ as $n\to\infty$.
\end{prop}

\begin{pf}
Let $\mu_i=\mu_i(n)$ be the number of $i\times i$ symmetric incidence matrices
with $n$ ones.
Let $s_k$ be the number of $k\times k$ symmetric matrices with $n$ ones,
given by
\begin{equation}\label{sdef}
s_k=\sum_{j=0}^{\lfloor n/2\rfloor}{{k\choose2}\choose j}{k\choose n-2j},
\end{equation}
where $j$ represents the number of ones off of the diagonal.
Then 
\[
s_k=\sum_{i=1}^k{k\choose i}\mu_i
\]
and
\[
S_{11}(n)=\sum_{i=1}^n\mu_i=\sum_{k=1}^\infty\frac{s_k}{2^{k+1}}
\]
by (\ref{sum1}), leading to
\[
S_{11}(n)
=\sum_{k=0}^\infty\frac1{2^{k+1}}\sum_{j=0}^\infty{{k\choose 2}\choose
j}{k\choose n-2j}\;.
\]
To compute this sum asymptotically, we approximate for $m\gg l\gg 1$
\[
{m\choose l}\sim\frac{m^l}{l!}\ee^{-l^2/2m}\;.
\]
The sums are dominated near $k\approx n/\log 2$ and $j\approx(n-\sqrt n)/2$, so
that we can justify
replacing the binomial coefficients by this approximation. We get
\begin{eqnarray*}
S_{11}(n)&\sim&
\sum_{k=0}^\infty\frac1{2^{k+1}}\sum_{j=0}^\infty\frac{{k\choose 2}^j}{j!}\frac{k^{n-2j}}{(n-2j)!}
\ee^{-j^2/2{k\choose 2}-(n-2j)^2/2k}\\
&\sim&
\sum_{k=0}^\infty\frac{k^n}{2^{k+1}}\sum_{j=0}^\infty\frac1{j!(n-2j)!2^j}\ee^{-j/k-j^2/k^2-(n-2j)^2/2k}\;,
\end{eqnarray*}
where in the last step we also replaced $(1-1/k)^j\sim \ee^{-j/k}$. 
Due to the concentration of the sum near $k\approx n/\log 2$ and 
$j\approx(n-\sqrt n)/2$, the argument of the exponential can be replaced by
\[
\ee^{-j/k-j^2/k^2-(n-2j)^2/2k}\sim \ee^{-\log2-(\log 2/2)^2}
={\textstyle{\frac{1}{2}}}\ee^{-(\log2)^2/4}\;=C_s,
\]
where $C_s$ is as in the Proposition. Identifying
\[
P(n)=\sum_{k=0}^\infty\frac{k^n}{2^{k+1}}
\]
and
\[
I(n)=\sum_{j=0}^\infty\frac{n!}{j!(n-2j)!2^j}\;,
\]
we arrive at
\[
S_{11}(n)\sim C_s\cdot I(n)P(n)/n!\;.
\]
\qed
\end{pf}

One may also generate matrices from $S_{11}$ uniformly at random.

Define an integer valued probability distribution function
$\psi_k$ by
\[
\psi_k=\frac{s_k2^{-k-1}}{S_{11}(n)}.
\]
\begin{theorem}\label{Sgen}
The following algorithm generates a random incidence matrix counted by
$S_{11}(n)$.
\begin{enumerate}
\item[(i)] Draw integer-valued random variables $K$
according to the probability distribution $\psi_k$.

\item[(ii)] Choose a $K\times K$ symmetric zero-one matrix
with $n$ ones and $K^2-n$ zeros uniformly at random.

\item[(iii)] Delete all rows and columns for which all entries are zero.
\end{enumerate}
\end{theorem}
The proof of Theorem~\ref{Sgen} is similar to the proof
of Theorem~\ref{Fgen}.

\bigskip

In general, we have $\Phi_{ij}(n)=\frac{1}{2}(F_{ijij}(n)+S_{ij}(n))$, where
\begin{itemize}
\item if $i=1$, then $S_{ij}(n)$ is the number of symmetric matrices with
$n$ ones and no zero rows, where repeated rows are forbidden if $j=0$ and
permitted if $j=1$;
\item if $i=0$, then $S_{ij}(n)$ is the number of classes of matrices
with $n$ ones and no zero rows (up to row and column permutations) which
are closed under transposition, with the same interpretation of $j$ as in
the other case.
\end{itemize}

We do not yet have asymptotics for these. It seems likely that, in all four
cases, $S_{ij}(n)=o(F_{ijij}(n))$, so that
$\Phi_{ij}(n)\sim\frac{1}{2}F_{ijij}(n)$. Table~\ref{sfns} gives some values
of these functions.

\begin{table}[htbp]
\caption{\label{sfns}Some counts for symmetric matrices and classes}
\[
\begin{array}{||c||r|r|r|r|r|r|r|r|r|r||}
\hline
n & 1 & 2 & 3 & 4 & 5 & 6 & 7 & 8 & 9 & 10\\
\hline
S_{00}(n) & 1 & 1 & 2 &  2 &  3 &   6 &&&& \\
S_{01}(n) & 1 & 1 & 2 &  4 &  6 &  10 &&&& \\
S_{10}(n) & 1 & 2 & 6 & 16 & 56 & 214 &  866 & 3796 & 17468 & \\
S_{11}(n) & 1 & 2 & 6 & 20 & 74 & 302 & 1314 & 6122 & 29982 & 154718 \\
\hline
\end{array}
\]
\end{table}

\section{The function $F_{1011}(n)$}

Recall that the number of incidence matrices with $n$ ones, no repeated
rows and matrices equal by row or column permutations unidentified
is denoted by $F_{1011}(n)$.
In this section we will show
\begin{theorem}
We have
$$
F_{1011}(n)=o\left(F_{1111}(n)\right).
$$
\end{theorem}
\begin{pf}
We will use the probabilistic method and the notation
used in the proof of 
Theorem~\ref{F1111}. The idea behind the proof is to show that
the probability tends to $0$
 that a randomly chosen incidence matrix counted by
$F_{1111}(n)$ does not have two rows with each containing all zeroes
except for a single one in the same column.

Define $E_{i,j}$, $1\leq i,j\leq n$, to be the event that both
$\{i\}$ and $\{j\}$ are blocks in the first preorder and that $i$ and
$j$ belong to the same block of the second preorder. When
$W=0$, $E_{i,j}$ corresponds to the event that the rows corresponding
to the blocks containing $i$ and $j$ in the incidence matrix are
different and contain unique ones appearing in the same column.

Let $P(n,k)$ be the number of preorders on $n$ elements with $k$ blocks.
Given a power series $f(z)=\sum_{n=0}^\infty f_n z^n$, define 
$[z^n]f(z)=f_n$.
We find that for any $1\leq i<j\leq n$,
\begin{eqnarray}\label{probE}
\mathbb{P}(E_{i,j})&=&
\sum_{k=1}^{n-2}\frac{P(n-2,k)}{P(n)}(k+2)(k+1)\cdot\frac{P(n-1)}{P(n)}
\nonumber\\
&=&
\frac{P(n-1)}{P(n)^2}\sum_{k=1}^{n-2}
(k+2)(k+1)P(n-2,k)
\end{eqnarray}

Using Lemma~1.1 of \cite{cs}, we find that
\begin{eqnarray*}
\sum_{k=1}^{n-2}
(k+2)(k+1)P(n-2,k)
&=&
(n-2)!
[z^{n-2}]\left(\sum_{n=0}^\infty(n+2)(n+1)(\ee^z-1)^n\right)\\
&=&
(n-2)!
[z^{n-2}]
\frac{d^2}{du^2}
\left(
\frac{u^2}{1-u}\,\Big|_{u=\ee^z-1}
\right).
\end{eqnarray*}
When singularity analysis (see Section 11 of \cite{aem}) can be applied,
as in this case,
the asymptotics of the coefficients of a generating function are determined
by the degree of its pole of smallest modulus. Therefore,
$$
\sum_{k=1}^{n-2}
(k+2)(k+1)P(n-2,k)
\sim
(n-2)!
[z^{n-2}]
\left(
\frac{2(\ee^z-1)^2}{(2-\ee^z)^3}
\right).
$$
The singularity of smallest modulus of $\left(2-\ee^z\right)^{-1}$
occurs at $z=\log 2$ with residue
$$
\lim_{z\to\log2}\left(\frac{z-\log2}{2-\ee^z}\right)
=\lim_{z\to\log2}\left(\frac{1}{-\ee^z}\right)=-\frac{1}{2},
$$
by l'H\^opital's rule. Hence,
$$
\sum_{k=1}^{n-2}
(k+2)(k+1)P(n-2,k)
\sim
\frac{(n-2)!}{4}
[z^{n-2}](\log 2-z)^{-3}
$$
from which singularity analysis and (\ref{Pasymp}) give
\begin{equation}\label{coeffasymp}
\sum_{k=1}^{n-2}
(k+2)(k+1)P(n-2,k)
\sim
\frac{(n-2)!}{8}
(\log 2)^{-n-1}n^2.
\end{equation}
The result of using (\ref{coeffasymp}) in (\ref{probE}) is
\begin{equation}\label{sim1}
\mathbb{P}(E_{i,j})\sim\frac{\log 2}{4n}.
\end{equation}
Define $X$ to be
$$
X=\sum_{1\leq i<j\leq n} I_{E_{i,j}},
$$
so that, conditional on the event $\{W=0\}$,
the event $\{X>0\}$ implies that the incidence matrix
produced by the algorithm has repeated rows.
The expectation of $X$ is
\begin{equation}\label{expX}
\mathbb{E}(X)={n\choose 2}\mathbb{P}(E_{i,j})\sim \frac{\log 2}{8} n
\end{equation}

We will next show that 
\begin{equation}\label{sim2}
\mathbb{E}\left((W)_r\cap E_{i,j}\right)
\sim\left(\frac{(\log 2)^2}{2}\right)^r\frac{\log 2}{4n}.
\end{equation}
The analog of (\ref{falling}) is
\begin{eqnarray*}
\mathbb{E}\left((W)_r\cap E_{i,j}\right)
&=&
\mathbb{E}\left(
\sum_{
\stackrel
{{\rm pairs \ }(i_s,j_s){\rm \ different}}
{i_s\neq i, \ j_s\neq j}
}
I_{i_1,j_1}\cdots I_{i_r,j_r}I_{E_{i,j}}\right)\\
&=&
\mathbb{E}\left(\sum_{
\stackrel
{{\rm all \ }i_s {\rm \ and \ } j_s{\rm \ different}}
{i_s\neq i, \ j_s\neq j}
}
I_{i_1,j_1}\cdots I_{i_r,j_r}I_{E_{i,j}}\right)\\
&&+\,
\mathbb{E}\left(\sum\nolimits^{\ast\ast}
I_{i_1,j_1}\cdots I_{i_r,j_r}I_{E_{i,j}}\right)
\end{eqnarray*}
with ${\displaystyle \sum\nolimits^{\ast\ast}}$ defined
to be the sum with all pairs $(i_s,j_s)$ different, but not all
$i_s, j_s$ different.

The first term corresponds to two preorders formed in the following
way. The $i_s,j_s$ are first selected.
One preorder is formed from the set of elements
$[n]\setminus\{i,j,j_1,j_2,\ldots,j_r\}$, the element
$j_s$ is added to the block containing $i_s$ for each $s$,
and then blocks $\{i\}$ and $\{j\}$ are inserted in the preorder.
Another preorder is formed from the set of elements
$[n]\setminus\{j,j_1,j_2,\ldots,j_r\}$, 
the element
$j_s$ is added to the block containing $i_s$ for each $s$, 
and then the element $j$ is added to the block containing $i$.
As a result,
\begin{eqnarray*}
&&
\mathbb{E}\left(
\sum_{
\stackrel
{{\rm pairs \ }(i_s,j_s){\rm \ different}}
{i_s\neq i, \ j_s\neq j}
}
I_{i_1,j_1}\cdots I_{i_r,j_r}\right)\\
&=&
\frac{(n-2)!}{2^r(n-2-2r)!}
\sum_{k=1}^{n-2-r}\frac{P(n-2-r,k)}{P(n)}(k+2)(k+1)\frac{P(n-1-r)}{P(n)}\\
&\sim&
\left(\frac{(\log 2)^2}{2}\right)^r\frac{\log 2}{4n}
\end{eqnarray*}
where we have used (\ref{Pasymp}) and the asymptotic form (\ref{coeffasymp}).

The second term is bounded using the same method that was used
to bound ${\displaystyle \sum\nolimits^\ast}$ 
in the proof of Theorem~\ref{F1111}. 
Letting $G^{\prime\prime}$ be an index over graphs on $n$ vertices
with two labelled disconnected
vertices $i$ and $j$ and $n-2$
unlabelled vertices which is not a matching on the unlabelled vertices, 
we have
\begin{eqnarray*}
&&
\mathbb{E}\left(\sum\nolimits^{\ast\ast}I_{i_1,j_1}\cdots I_{i_r,j_r}I_{E_{i,j}}\right)\\
&\leq&
\sum_{G^{\prime\prime}} (n-2)^v \sum_k\frac{P(n-2-v+c,k)}{P(n)}(k+1)(k+2)
\frac{P(n-1-v+c)}{P(n)}\\
&=&\sum_{G^{\prime\prime}}O\left(n^{2c-v-1}\right)\\
&=&O(n^{-2}).
\end{eqnarray*}
Consequently we have shown (\ref{sim2}).

The asymptotics (\ref{sim1}) and (\ref{sim2}) and method of moments
argument giving (\ref{W0}) imply that
\[
\mathbb{P}(W=0\mid E_{i,j})\sim \exp\left(-\frac{(\log 2)^2}{2}\right)
\]
and therefore an application of Bayes' Theorem with (\ref{W0})
and (\ref{sim1}) results in
\[
\mathbb{P}(E_{i,j}\mid W=0)\sim\frac{\log 2}{4n}.
\]
The observations above result in
\begin{equation}\label{Xexpcond}
\mathbb{E}(X\mid W=0)\sim \frac{\log 2}{8} n.
\end{equation}
Comparison of (\ref{expX}) and (\ref{Xexpcond}) makes it clear that
that conditioning on the event
$\{W=0\}$ does not asymptotically
affect the expectation of $X$.

In a similar way we can find the asymptotics of the conditional second
falling moment $\mathbb{E}(X(X-1)\mid Y=0)$. The unconditioned second moment equals
\begin{eqnarray*}
\mathbb{E}(X(X-1))&=&\sum_{(i_1,j_1)\neq (i_2,j_2)} \mathbb{P}(E_{i_1,j_1}\cap E_{i_2,j_2})\\
&=&
\sum_{\{i_1,j_1\}\cap \{i_2,j_2\}=\emptyset} \mathbb{P}(E_{i_1,j_1}\cap E_{i_2,j_2})+
\sum_{|\{i_1,j_1\}\cap \{i_2,j_2\}|=1} \mathbb{P}(E_{i_1,j_1}\cap E_{i_2,j_2})\\
&=&
\frac{n!}{4(n-4)!} \sum_{k=1}^{n-4} \frac{P(n-4,k)}{P(n)}(k+1)(k+2)(k+3)(k+4)
\frac{P(n-2)}{P(n)}\\
&&+\,O\left(n^3\sum_{k=1}^{n-4} \frac{P(n-3,k)}{P(n)}(k+1)(k+2)(k+3)
\frac{P(n-2)}{P(n)}\right)
\end{eqnarray*}
An application of singularity analysis as used to derive (\ref{coeffasymp})
produces
\[
\mathbb{E}(X(X-1))=\frac{(\log 2)^2}{64}n^2+O(n).
\]
Arguing as we did for $\mathbb{E}(X\mid W=0)$ shows that
\[
\mathbb{E}(X(X-1)\mid W=0)=\frac{(\log 2)^2}{64}n^2+O(n);
\]
we omit the details.

The variance of $X$ conditioned on $W=0$ is
\begin{eqnarray}\label{varXcond}
{\rm Var}(X\mid W=0)&=&\mathbb{E}(X(X-1)\mid W=0)+\mathbb{E}(X\mid W=0)
-(\mathbb{E}(X\mid W=0))^2\nonumber\\
&=& o(n^2).
\end{eqnarray}
Chebyshev's inequality applied with (\ref{Xexpcond}) and (\ref{varXcond})
now gives
\[
\mathbb{P}(X=0\mid W=0)=o(1).
\]
Hence, an asymptotically insignificant fraction of incidence
matrices do not have repeated rows which implies that $F_{1011}(n)
=o\left(F_{1111}(n)\right)$.
\qed
\end{pf}

\section{The functions $F_{0011}(n)$ and $F_{0111}(n)$}

The function $F_{0111}(n)$ counts vertex-labelled hypergraphs on $n$ vertices,
while $F_{0011}(n)$ counts the simple vertex-labelled hypergraphs. For
completeness, we include the formulae from the work of Martin 
Klazar~\cite{klazar}.

\begin{theorem}
\begin{description}
\item{(a)} For all $n$, we have
\begin{eqnarray*}
F_{0011}(n)&=&\sum_{\lambda\vdash n}\sum_{j=l}^n
\prod_{i=1}^l{{j\choose i}\choose a_i}\sum_{m=j}^n(-1)^{m-j}{m\choose j},\\
F_{0111}(n)&=&\sum_{\lambda\vdash n}\sum_{j=l}^n
\prod_{i=1}^l{{j\choose i}+a_i-1\choose a_i}\sum_{m=j}^n(-1)^{m-j}{m\choose j},
\end{eqnarray*}
where $\lambda=1^{a_1}2^{a_2}\cdots l^{a_l}$ is a partition of $n$ with 
$a_l>0$.
\item{(b)} For all $n$, we have
\[F_{0011}(n)\le F_{0111}(n)\le 2F_{0011}(n).\qquad\Box\]
\end{description}
\end{theorem}

Part (b) raises the question of whether $F_{0011}(n)/F_{0111}(n)$ tends to
a limit as $n\to\infty$, and particular, whether the limit is~$1$ (that is,
whether almost all labelled hypergraphs are simple).

The paper~\cite{klazar} also gives recurrence relations for the two functions.
Klazar subsequently showed~\cite{klazar2} that both functions are 
asymptotically
\[(1/\log 2+o(1))^nb(n),\]
where $b(n)$ is the $n$th Bell number (the number of partitions of 
$\{1,\ldots,n\}$). Details of the asymptotics of $b(n)$ can be found
in~\cite{aem}. In particular, since $b(n)/n!=(1+o(1))^n$, we see that
$F_{0111}(n)=o(F_{1111}(n))$, and in fact
\[F_{0111}(n) = (\log 2+o(1))^n F_{1111}(n).\]
This and the result of the last section, together with the facts that
$F_{1101}(n)=F_{0111}(n)$ and $F_{1110}(n)=F_{1011}(n)$ and that $F_{ijkl}(n)$
is monotone increasing in each of $i,j,k,l$, justify our earlier claim that
$F_{ijkl}(n)=o(F_{1111}(n))$ for $(i,j,k,l)\ne(1,1,1,1)$.

\section{A rough lower bound for $F_{0001}(n)$}

The number $F_{0101}(n)$ of unlabelled hypergraphs with weight~$n$ is
not smaller than the number of graphs with $n/2$ edges and no isolated
vertices. We show that this number grows faster than exponentially. In
fact, our argument applies to $F_{0001}$, since we use simple graphs.

Consider simple graphs with $m$ vertices and $n$ edges, where $m=o(n)$ and
$n=o(m^2)$. The number of such graphs, up to isomorphism, is at least
\[\frac{{m(m-1)/2\choose n}}{m!} > \frac{(cm^2)^n}{n!\,m!}
> \frac{c^nm^{2n}}{n^nm^m}\]
for some constant $c$.
Put $F$ equal to the logarithm of the right-hand side:
\[F=c'n+2n\log m-n\log n-m\log m\]
for some constant $c'$. 
Putting $m=c''n/\log n$, for some constant $c''$, we get
\[F = n\log n - 2n\log\log n + O(n).\]

We conclude:
\begin{prop}
For any $\epsilon>0$, we have
\[F_{0001}(n) \ge \left(\frac{n}{(\log n)^{2+\epsilon}}\right)^n\]
for $n\ge n_0(\epsilon)$.
\end{prop}

\paragraph{Remark} The asymptotics of the number of graphs with
no isolated vertices, having a given number of vertices and edges,
has a long history: see~\cite{wright} for an early paper on this topic, 
and \cite{BCM} for a recent result.


\begin{thebibliography}{99}

\bibitem{Bar}
J. P. Barthelemy, An asymptotic equivalent for the number
of total preorders on a finite set, 
\textit{Discrete Math.} \textbf{29} (1980) 311--313.

\bibitem{BCM}
E. A. Bender, E. R. Canfield and B. D. McKay,
The asymptotic number of labeled graphs with $n$ vertices, $q$ edges, and 
no isolated vertices,
\textit{J. Combinatorial Theory} (A) \textbf{80} (1997), 124--150.

\bibitem{Ben}
E. A. Bender and S. G. Williamson,
\textit{Foundations of Combinatorics with Applications},\hfil\break
\texttt{http://math.ucsd.edu/\lower0.5ex\hbox{\~{}}ebender/CombText/}

\bibitem{C}
P. J. Cameron, 
\textit{Combinatorics: Topics, Techniques, Algorithms},
Cambridge University Press, 1994.

\bibitem{cgm}
P. J. Cameron, D. A. Gewurz and F. Merola, Product action,
to appear.

\bibitem{cps}
P. J. Cameron, Thomas Prellberg and Dudley Stark, 
Asymptotic enumeration of incidence matrices,
\textit{Journal of Physics (Conference Series)}, to appear.

\bibitem{cs}
P. J. Cameron and Dudley Stark, Random preorders, submitted.

\bibitem{gap}
The GAP Group,
\textsf{GAP} --- Groups, Algorithms, and Programming, Version 4.3, 2002,
\texttt{http://www.gap-system.org}

\bibitem{klazar}
M. Klazar,
Extremal problems for ordered hypergraphs: small patterns and some enumeration,
\textit{Discrete Appl. Math.} \textbf{143} (2004), 144--154.

\bibitem{klazar2}
M. Klazar,
Counting set systems by weight,
\textit{Electr. J. Combinatorics} \textbf{12(1)} (2005), \#R11, (8pp).

\bibitem{MB}
H. Maassen and T. Bezembinder,
Generating random weak orders and the probability of a Condorcet winner.
\textit{Soc. Choice Welfare} \textbf{19} (2002) 517--532.

\bibitem{MM}
M. Maia and M. Mendez,
On the arithmetic product of combinatorial species,
preprint available from
\texttt{arXiv:math.CO/0503436}

\bibitem{aem}
A. M. Odlyzko,
Asymptotic enumeration methods,
In R.~L.~Graham, M.~Gr\"otschel and L.~Lov\'asz (eds.), 
\textit{Handbook of Combinatorics}, Vol. 2, North-Holland, Amsterdam, 1995,
pp.~1063--1229.

\bibitem{oeis}
N. J. A. Sloane (ed.), \textit{The On-Line Encyclopedia of Integer 
Sequences},\hfil\break
\texttt{http://www.research.att.com/\lower0.5ex\hbox{\~{}}njas/sequences/}

\bibitem{wright}
E. M. Wright,
Two problems in the enumeration of unlabelled graphs.
\textit{Discrete Math.} \textbf{9} (1974), 289--292.

\end{thebibliography}
\end{document}